\documentclass[12pt]{amsart}
\usepackage{amssymb,amsthm}
\usepackage[all]{xy}

\numberwithin{equation}{section}
\theoremstyle{plain}
\newtheorem{thm}{Theorem}[section]

\newtheorem{prop}[equation]{Proposition}

\newtheorem{questions}[equation]{Questions}

\theoremstyle{remark}

\theoremstyle{definition}

\newtheorem{ex}[equation]{Example}

\DeclareMathOperator{\HH}{HH}
\DeclareMathOperator{\chara}{char}
\DeclareMathOperator{\coh}{H}
\DeclareMathOperator{\cx}{cx}
\DeclareMathOperator{\Hom}{Hom}
\DeclareMathOperator{\Ext}{Ext}

\newcommand{\C}{\mathbb{C}}
\newcommand{\Z}{\mathbb{Z}}

\newcommand{\ot}{\otimes}
\newcommand{\DOT}{\setlength{\unitlength}{1pt}\begin{picture}
              (2.5,2)(1,1)\put(2,3.5){\circle*{2}}\end{picture}}
\newcommand{\bu}{\DOT}

\title[Varieties for modules]
{Varieties for modules\\ of finite dimensional Hopf algebras}
\author{Sarah Witherspoon}
\address{Department of Mathematics, Texas A\&M University, College Station, Texas 77843, USA}
\email{sjw@math.tamu.edu}
\date{2 December 2016}

\thanks{The author was partially supported by
NSF grant \#DMS-1401016.}

\dedicatory{Dedicated to Professor David J.\ Benson on the
occasion of his 60th birthday.}

\begin{document}

\begin{abstract}
We survey variety theory for modules of finite dimensional Hopf algebras, recalling
some definitions and basic properties of 
support and rank varieties where they are known. 
We focus specifically on properties known for
classes of examples such as finite group algebras and finite group schemes. 
We list open questions about 
tensor products of modules and projectivity, where varieties may play a role
in finding answers. 
\end{abstract}

\maketitle

\section{Introduction}

For a given group or ring, one wants to understand its representations in a meaningful way.
It is often too much to ask for a full classification of all indecomposable modules,
since one may work in a setting of wild representation type.
Varieties can then be an important tool for organizing 
representations and extracting information.
In the theory of varieties for modules, one associates to each
module a geometric space---typically an affine or projective variety---in such
a way that representation theoretic properties are encoded in the space.  
Varieties for modules originated in finite group representation 
theory, in work of Quillen~\cite{Quillen} and Carlson~\cite{Carlson}.
This theory and all required background material 
is elegantly presented in Benson's book~\cite{benson91b}. 
The theory has been adapted to many other settings, such as finite group schemes, 
algebraic groups, Lie superalgebras, quantum groups, and self-injective algebras.
See, e.g., 
\cite{AB,AI,BKN,DS,EHSST,Friedlander,FP86,Jantzen,Ostrik,PS,SFB,Snashall-Solberg}.

In this survey article, we focus on finite dimensional Hopf algebras,
exploring the boundary between those
whose  variety theory behaves as one expects, arriving from finite
group representation theory, and those where it does not.
We give definitions of support varieties in terms of Hochschild cohomology
from~\cite{EHSST,Snashall-Solberg}, and in terms of
Hopf algebra cohomology as a direct generalization of group 
cohomology from~\cite{Feldvoss-W,FP05,F-Pevtsova,Ostrik}. 
We recall which Hopf algebras are known to have finitely
generated cohomology, opening the door to these standard versions of support
varieties. 
We also briefly summarize the rank varieties which are defined representation theoretically in
a more limited array of settings, yet are indispensable where they
are defined. 

We are most interested in
the tensor product property, that is,
the property that the variety of a tensor product of modules is equal to the intersection
of their varieties. 
This is known to hold for modules of some Hopf algebras, known not to hold for others,
and is an open question for most. 
We  look at some related questions about
tensor products of modules: 
(i) If the tensor product of two modules in one order
is projective, what about their tensor product in the other order?
(ii) If a tensor power of a module is projective, need the module itself be projective?
The answers to both questions  are yes for finite group algebras and finite group schemes, 
while the answers to  both are no for some types of Hopf algebras, as we will see.
In fact,  any finite dimensional Hopf algebra
satisfying the tensor product property is a subalgebra of 
one that does not and for which 
the above two questions have negative answers.

The open questions we discuss in this article are important for gaining
a better understanding of the representation theory of finite dimensional Hopf algebras.
Their module categories enjoy a rich structure due to existence of tensor products. 
Varieties are a great tool for
understanding these tensor products when one understands the relationship between
them and their varieties.

Throughout, we will work over an algebraically closed field $k$,
although there are known results for more general fields and  ground rings in some contexts.
Sometimes we will assume $k$ has positive characteristic, and
sometimes that it has characteristic 0. 
All tensor products will be taken over $k$ unless otherwise indicated,
that is, $\ot = \ot_k$. 
All modules will  left modules, finite dimensional over $k$,
 unless otherwise stated. 

\section{Hopf algebras}

A {\em Hopf algebra} is an algebra $A$ over the field $k$
together with $k$-linear maps
$\Delta: A\rightarrow A\ot A$ (comultiplication), 
$\varepsilon: A\rightarrow k$ (counit or augmentation), and 
$S: A\rightarrow A$ (antipode or coinverse)  satisfying the following 
properties: 
The maps $\Delta$ and $\varepsilon$ are
algebra homomorphisms, and $S$ is an algebra anti-homomorphism
(i.e., it reverses the order of multiplication).
Symbolically writing $\Delta(a)=\sum a_1\ot a_2$ (Sweedler notation),
we also require 
$$
\begin{aligned}
   (1\ot \Delta) (\Delta(a)) & =  (\Delta\ot 1)(\Delta(a))
   \ \ \mbox{ (coassociativity)}, \\
\sum \varepsilon (a_1) a_2  = & \ a = \sum a_1\varepsilon(a_2)
   \ \ \mbox{ (counit property)}, \\
\sum S(a_1)a_2 = & \ \varepsilon(a)\cdot 1 =\sum a_1 S(a_2)
  \ \ \mbox{ (antipode property)} 
\end{aligned}
$$
for all $a\in A$. 
We say that $A$ is {\em cocommutative} if $\tau\circ \Delta = \Delta$, where 
$\tau:A\ot A\rightarrow A\ot A$ is the twist map, that is, $\tau(a\ot b) = b\ot a$
for all $a,b\in A$. 
For more details, see, e.g., \cite{montgomery93}. 

Standard examples of Hopf algebras, some of which will 
reappear in later sections, are: 

\begin{ex}
$A=kG$, the group algebra of a finite group $G$, with $\Delta(g)=g\ot g$, $\varepsilon(g)=1$,
and $S(g)=g^{-1}$ for all $g\in G$. This Hopf algebra is cocommutative.
\end{ex}

\begin{ex}\label{ex:dual}
$A= k[G]=\Hom_k(kG,k)$, the linear dual of the group algebra $kG$,
in which multiplication is pointwise on group
elements, that is, $(ff')(g) = f(g) f'(g)$ for all $g\in G$ and $f,f'\in k[G]$. 
Comultiplication is given as follows.
Let $\{p_g\mid g\in G\}$ be the basis of $k[G]$ dual to $G$.
Then 
   $$\Delta(p_g)=\sum_{\substack{a,b\in G \\    ab=g}} p_a\ot p_b,$$
$\varepsilon(p_g)=\delta_{g,1}$, and $S(p_g)=p_{g^{-1}}$ for all $g\in G$.
This Hopf algebra is noncocommutative when $G$ is nonabelian.  
\end{ex}

\begin{ex}
 $A=U({\mathfrak{g}})$, the universal enveloping algebra of
a Lie algebra $\mathfrak g$, with 
$\Delta(x)=x\ot 1 + 1\ot x$, $\varepsilon(x)=0$, and $S(x) = -x$ 
for all $x\in {\mathfrak{g}}$. The maps $\Delta$ and $\varepsilon$
are extended to be algebra homomorphisms, and $S$ to be an algebra
anti-homomorphism. This is an infinite dimensional cocommutative Hopf algebra.
In case the characteristic of $k$ is a prime $p$, and $\mathfrak g$
is a restricted Lie algebra, its restricted enveloping algebra
$u({\mathfrak {g}})$ is a finite dimensional cocommutative Hopf algebra
with analogous comultiplication, counit, and antipode.  
\end{ex}

\begin{ex}\label{ex:sqg}
 $A = U_q({\mathfrak{g}})$ or $A=u_q({\mathfrak{g}})$,
the  infinite dimensional 
quantum enveloping algebras and some finite dimensional versions
(the small quantum groups).
See, e.g.,~\cite{ginzburg-kumar93} for the definition in the general case.
Here we  give just one small example explicitly: Let $q$ be a primitive
complex $n$th root of unity, $n>2$. Let
$u_q({\mathfrak{sl}}_2)$ be the $\mathbb{C}$-algebra generated by
$E,F,K$ with $E^n=0$, $F^n=0$, $K^n=1$, $KE=q^2 EK$,
$KF=q^{-2}FK$, and $$EF=FE + \frac{K-K^{-1}}{q-q^{-1}}.$$ 
Let $\Delta(E) = E\ot 1 + K\ot E$, $\Delta(F) = F\ot K^{-1} + 1\ot F$,
$\Delta(K) =K\ot K$, $\varepsilon(E)=0$, $\varepsilon(F) =0$,
$\varepsilon(K) = 1$, $S(E) = -K^{-1}E$, $S(F) = -FK$, and $S(K)=K^{-1}$. 
This is a finite dimensional noncocommutative Hopf algebra. 
\end{ex}

\begin{ex}\label{ex:QEAG}
$A$ is a  quantum elementary abelian group: 
Let $m$ and $n$ be positive integers, $n\geq 2$. 
Let $q$ be a primitive complex $n$th root of unity, and let $A$ be the $\C$-algebra
generated by $x_1,\ldots, x_m, g_1,\ldots, g_m$ with relations 
$x_i^n=0$, $g_i^n=1$, $x_ix_j=x_jx_i$, $g_ig_j=g_jg_i$, 
and $g_ix_j = q^{\delta_{i,j}}  x_j g_i$ for all $i,j$.
Comultiplication is given by $\Delta(x_i)= x_i\ot 1 + g_i\ot x_i$,
$\Delta(g_i)=g_i\ot g_i$, counit $\varepsilon(x_i)=0$, $\varepsilon(g_i)=1$,
and antipode $S(x_i)= - g_i^{-1}x_i$, $S(g_i)= g_i^{-1}$ for all $i$.
This Hopf algebra is finite dimensional and  noncocommutative. 
\end{ex}

We return to the general setting of a Hopf algebra $A$, which we 
assume from now on is finite dimensional over $k$. 
Letting $M$ and $N$ be $A$-modules, their tensor product $M\ot N$ is
again an $A$-module via the comultiplication map $\Delta$, that is, 
 $a\cdot (m\ot n ) = \sum (a_1\cdot m) \ot (a_2\cdot n)$
for all $a\in A$, $m\in M$, $n\in N$. 
The category of (finite dimensional) 
$A$-modules is  a rigid tensor category: 
There is a unit object given by the field $k$ under action via
the counit $\varepsilon$, i.e., 
$a\cdot c = \varepsilon(a)c$ for all $a\in A$ and $c\in k$.
There are dual objects: Let $M$ be a finite dimensional $A$-module, and
let $M^*= \Hom_k (M,k)$, an $A$-module via $S$: 
$(a\cdot f)(m) = f( S(a)\cdot m)$ for all $a\in A$, $m\in M$.
See, e.g.,~\cite{BK} for details on rigid tensor categories. 

The following proposition is proven  in~\cite[Proposition 3.1.5]{benson91a}.
An alternative proof is to observe  that 
$\Hom_A(P,\Hom_k(M, - )) \cong \Hom_A(P\otimes M, - )$
as functors, where the action of $A$ on $\Hom_k (M,N)$, for  $A$-modules $M,N$,
is given by $(a\cdot f)(m)
= \sum a_1 \cdot ( f(S(a_2)\cdot m))$ for all $a\in A$, $m\in M$, and $f\in \Hom_k(M,N)$. 
(A similar argument applies to $M\ot P$.) 

\begin{prop}\label{projective}
If $P$ is a projective  $A$-module, and $M$ is any $A$-module, then
both $P\ot M$ and $M\ot P$ are  projective  $A$-modules. 
\end{prop}

In Section~\ref{sec:questions}, we will consider other connections between
projectivity and tensor products. 

\quad

\section{Varieties for modules}\label{sec:var}

There are many versions of varieties for modules, depending on which rings and 
modules are of interest.
Here we will present the support variety theory of Erdmann, Holloway, Snashall, Solberg,
and Taillefer~\cite{EHSST,Snashall-Solberg} for self-injective algebras
(based on Hochschild cohomology), as well as the closely related  
generalization to Hopf algebras of support varieties for finite group 
representations (see, e.g.,~\cite{Feldvoss-W,FP,FP05,Ostrik}). 
See also Solberg's excellent survey~\cite{Solberg} for more details. 

\subsection*{Hochschild cohomology}
Let $A$ be an associative $k$-algebra. 
Let $A^e= A\ot A^{op}$, with $A^{op}$ the opposite algebra to $A$.
Consider $A$ to be an $A^e$-module via left
and right multiplication, that is, $(a\ot b)\cdot c = acb$ for all $a,b,c\in A$.
The {\em Hochschild cohomology} of $A$ is
$$
   \HH^*(A)= \Ext^*_{A^e}(A,A).
$$
The graded vector space $\HH^*(A)$ is a graded commutative ring under
Yoneda composition/cup product~\cite{benson91b,Suarez-Alvarez}, and
 $\HH^0(A)\cong Z(A)$, the center of $A$.   
If $M$ is an $A$-module, then the Hochschild cohomology ring $\HH^*(A)$
acts on $\Ext^*_A(M,M)$ via $-\ot_A M$ followed by Yoneda composition.

Now suppose $A$ is finite dimensional and self-injective.
For example, a finite dimensional Hopf algebra is a
Frobenius algebra, and therefore is self-injective. 
We will make  some assumptions,  as in \cite{EHSST}:

\quad

\noindent
Assume there is a graded subalgebra $H$ of $\HH^*(A)$ such that
\begin{itemize}
\item[(fg1)]  $ H$  is finitely generated, commutative, and $H^0=\HH^0(A)$, and 
\item[(fg2)] for all finite dimensional $A$-modules $M$, the Ext space 
$\Ext^*_A(M,M)$ is  finitely generated as an $H$-module.
\end{itemize}

\smallskip

For a finite dimensional $A$-module $M$, let
$I_A(M)$ be the annihilator in $H$ of $\Ext^*_A(M,M)$. The {\em support variety} of $M$ is 
\begin{equation}\label{eqn:suppvar}
    V_A(M) = {\rm{Max}} ( H/I_A(M)),
\end{equation}
the maximal ideal spectrum of $H/I_A(M)$.
This is the set of maximal ideals as a topological space under the Zariski
topology.
Alternatively, one considers homogeneous prime ideals as in
some of the given references.
Here we choose  maximal ideals as in~\cite{EHSST}.

\begin{ex}\label{ex:gpalg}
Let $A=kG$, where $G$ is a finite group.
We assume the characteristic of $k$ is a prime $p$ dividing the
order of $G$, since otherwise $kG$ is semisimple by Maschke's Theorem. 
The {\em group cohomology} ring is $\coh^*(G,k)= \Ext^*_{kG}(k,k)$.
More generally if $M$ is a $kG$-module, set
$\coh^*(G,M)= \Ext^*_{kG}(k,M)$.
There is an  algebra isomorphism:
$$
   \HH^*(kG)\cong \coh^*(G, (kG)^{ad})
$$
(see, e.g., \cite[Proposition~3.1]{siegel-witherspoon99}),
where the latter is group cohomology with coefficients in the adjoint $kG$-module $kG$
(on which $G$ acts by conjugation). 
 The group cohomology $\coh^*(G,k)$
then embeds into Hochschild cohomology $\HH^*(kG)$, 
since the trivial coefficients $k\cdot 1$ embed as a direct summand of $(kG)^{ad}$. 
Let $$H= \coh^{ev}(G,k)\cdot \HH^0(kG),$$ where $\coh^{ev}(G,k)$
is $\coh^*(G,k)$ if $\chara(k)=2$ and otherwise is
the subalgebra of $\coh^*(G,k)$ 
generated by its homogeneous even degree elements,
considered to be a subalgebra of Hochschild cohomology $\HH^*(kG)$
via the embedding discussed above. 
Then $H$ satisfies (fg1) and (fg2). 
The traditional definition of varieties for $kG$-modules uses
simply $\coh^{ev}(G,k)$ instead of $H$, the difference being 
the inclusion of the elements of $\HH^0(kG)\cong Z(kG)$.
If $G$ is a $p$-group, there is no difference in the theories
since $Z(kG)$ is local. 
See, e.g.,~\cite{benson91b} for details, 
including descriptions of the original work of
Golod~\cite{golod}, Venkov~\cite{venkov}, and Evens~\cite{Ev} on finite generation. 
If $G$ is not a $p$-group, the representation theoretic information contained in the
varieties will be largely the same in the two cases (the only exception being the
additional information of which block(s) a module lies in).
\end{ex}

Returning to the general setting of a finite dimensional self-injective algebra $A$, 
the support varieties defined above enjoy many useful properties \cite{EHSST}, some of which we collect
below. We will need to define the complexity of a module:
The {\em complexity} $\cx_A(M)$ of a finite dimensional
$A$-module $M$ is the rate of growth of a minimal
projective resolution. That is, if $P_{\bu}$ is a 
minimal projective resolution of $M$, then $\cx_A(M)$ is the
smallest nonnegative integer $c$ such that there is a real number $b$ and positive integer $m$ 
for which $\dim_k (P_n) \leq b n ^{c-1}$ for all  $n\geq m$. 
A projective module has complexity 0. 
The converse is also true, as stated in the following proposition.

\begin{prop}\label{prop:properties} \cite{EHSST,Snashall-Solberg}
Let $A$ be a finite dimensional self-injective algebra for
which there is an algebra $H$ satisfying (fg1) and (fg2).
Let $M$ and $N$ be finite dimensional $A$-modules. Then: 
\begin{itemize}
\item[(i)] $\dim V_A(M) = \cx_A(M)$.
\item[(ii)] $V_A(M\oplus N)=V_A(M)\cup V_A(N)$.  
\end{itemize}
Moreover, $\dim V_A(M) = 0$ if and only if $M$ is projective. 
\end{prop}

We will apply the support variety theory of~\cite{EHSST,Snashall-Solberg}, as outlined above,
to a finite dimensional Hopf algebra $A$, provided
there exists an algebra $H$ satisfying (fg1) and (fg2).

\subsection*{Hopf algebra cohomology}
Alternatively, one may generalize support varieties for finite groups
directly:
The {\em cohomology} of the Hopf algebra $A$ is $$\coh^*(A,k)= \Ext^*_A(k,k).$$ 
The cohomology $\coh^*(A,k)$ is a graded commutative ring under
Yoneda composition/cup product~\cite{Suarez-Alvarez}. 
If $M$ is an $A$-module,
consider $\Ext^*_A(M,M)$ to be an $\coh^*(A,k)$-module via $ - \ot M$ followed
by Yoneda composition. 
We make the following assumptions,
as in~\cite{Feldvoss-W}:

\medskip

\noindent
Assume that 
\begin{itemize}
\item[(fg1$'$)] $\coh^*(A,k)$ is a finitely generated algebra, and 
\item[(fg2$'$)] for all finite dimensional $A$-modules $M$, the Ext space 
$\Ext^*_A(M,M)$ is finitely generated as a module over 
$\coh^*(A,k)$.
\end{itemize}

\smallskip

Then one defines the {\em support variety} of an $A$-module $M$ to be 
the maximal ideal spectrum of the quotient of $\coh^*(A,k)$ by
the annihilator of $\Ext^*_A(M,M)$.
By abuse of notation, we will also write $V_A(M)$ for this variety,
and in the sequel it will be clear in each context  which is meant. 
If one wishes to work with a commutative ring from the beginning,
and not just a graded commutative ring,
then in characteristic not 2, 
one first restricts to the subalgebra $\coh^{ev}(A,k)$ of $\coh^*(A,k)$ generated
by all homogeneous elements of even degree. 
(The odd degree elements are nilpotent, and so the varieties are the same.) 
Proposition~\ref{prop:properties} holds for these varieties~\cite{Feldvoss-W}.

There is a close connection between this version of support variety and
that defined earlier via Hochschild cohomology:
Just as in Example~\ref{ex:gpalg}, Hopf algebra cohomology
$\coh^*(A,k)$ embeds into Hochschild cohomology $\HH^*(A)$.
See, e.g.,~\cite{ginzburg-kumar93} where this fact was first noted and
the appendix of~\cite{PW1} for a proof outline. 
One may then take $H= \coh^{ev}(A,k)\cdot \HH^0(A)$ in order to
define support varieties as in~(\ref{eqn:suppvar}). 
The only difference between these two versions of support variety
is the inclusion of the elements of $\HH^0(A)\cong Z(A)$.
Thus there is a finite surjective map from the variety defined
via the subalgebra $H$ of Hochschild cohomology to the variety 
defined via Hopf algebra cohomology. 

\subsection*{Rank varieties}
We now consider  the rank varieties that were first 
introduced by Carlson~\cite{Carlson} for studying finite group representations.
We recall his definition and discuss some Hopf algebras for which there are analogs.
Carlson's rank varieties are defined for elementary abelian $p$-groups.
For a finite group $G$, its elementary abelian $p$-subgroups
detect projectivity by Chouinard's Theorem~\cite[Theorem~5.2.4]{benson91b}, and
form the foundation for stratification of support varieties~\cite{AS82,Quillen}.
Thus it is important to understand the elementary abelian $p$-subgroups 
of $G$ and their rank varieties as defined below. 

Suppose $k$ is a field of prime characteristic $p$.
An {\em elementary abelian $p$-group} is a group of the form 
$E = (\Z/p\Z)^n$ for some $n$.
Write $E=\langle g_1,\ldots, g_n\rangle$, where $g_i$ generates
the $i$th copy of $\Z/p\Z$ in $E$.
For each $i$, let $x_i = g_i -1$, and note that $x_i^p=0$ since 
$\chara(k) =p$ and $g_i^p=1$.
It also follows that any element of the group algebra $kE$ of the form
$\lambda_1 x_1+ \cdots + \lambda_n x_n$ ($\lambda_i\in k$)
has $p$th power 0.
Thus for each choice of scalars $\lambda_1,\ldots,\lambda_n$,
 there is an algebra homomorphism
\begin{eqnarray*}
  k[t]/(t^p) & \rightarrow & kE \\
  t & \mapsto & \lambda_1 x_1 + \cdots + \lambda_n x_n  . 
\end{eqnarray*}
The image of this homomorphism is a subalgebra of $kE$ that
we will denote by $k \langle \lambda_1 x_1 + \cdots + \lambda_n x_n\rangle$.
Note that it is isomorphic to $k \Z/p\Z$ where the group $\Z/p\Z$ is
generated by $1+\lambda_1 x_1 + \cdots + \lambda_n x_n$.
The corresponding subgroup 
of the group algebra $kE$ is called a {\em cyclic shifted subgroup} of $E$. 
The {\em rank variety} of a $kE$-module $M$ is
\[
   V^r_E(M) = \{ {\mathbf{0}} \} \cup \{ (\lambda_1,\ldots, \lambda_n)\in k^n-\{ {\mathbf{0}}\}
   \mid M\!\downarrow_{k\langle \lambda_1x_1+\cdots +\lambda_nx_n\rangle}
    \mbox{ is not free} \} ,
\]
where the downarrow indicates restriction to the subalgebra.
Avrunin and Scott~\cite{AS82} proved that the rank variety $V^r_E(M)$
is homeomorphic to the support variety $V_E(M)$ (which we have
also denoted $V_{kE}(M)$).
 Information about a more general finite group $G$ is obtained by looking at
all its elementary abelian $p$-subgroups.
It is very useful to have on hand these rank  varieties for modules,
as another way to view the support varieties. 

Friedlander and Pevtsova~\cite{F-Pevtsova} generalized rank varieties to
finite dimensional cocommutative Hopf algebras $A$ (equivalently finite group
schemes), building on earlier work of Friedlander and Parshall~\cite{FP}
and Suslin, Friedlander, and Bendel~\cite{SFB}. 
The role of  cyclic shifted subgroups is played by subalgebras
isomorphic to
$k[t]/(t^p)$, or more generally by algebras $K[t]/(t^p)$ for field extensions
$K$ of $k$, and specific types of maps to extensions $A_K$.
A notion of rank variety for quantum elementary abelian groups is defined
in~\cite{PW1}, where the role of cyclic shifted subgroups is played by
subalgebras isomorphic to $k[t]/(t^n)$ with $n$  the order of
the root of unity $q$.
Scherotzke and Towers~\cite{SchTo} defined rank varieties for $u_q({\mathfrak{sl}}_2)$, and for the related Drinfeld doubles of Taft algebras, via certain subalgebras
detecting projectivity. 
Rank varieties have been defined as well for a number of algebras that 
are not Hopf algebras; see, e.g.,~\cite{BEH,BG,BE09}. 
In general though, it is not always clear what the right definition 
of rank variety should be, if any.

\section{Open questions and some positive answers}\label{sec:questions}

We next ask some questions about  finite
dimensional Hopf algebras, their representations, and varieties. 
We refer to the previous section for descriptions of the 
 support and rank varieties relevant to Question~\ref{q}(2) below.
We have purposely not specified choices of  varieties 
for the question, and  answers may depend on  choices. 
However, answers to the purely representation theoretic
Questions~\ref{q}(3) and (4) below do not. 

\begin{questions}\label{q}
Let $A$ be a finite dimensional Hopf algebra. 
\begin{itemize}
\item[(1)] Does $\coh^*(A,k)$ satisfy (fg1$\, '$) and (fg2$\, '$), or 
does there exist a subalgebra $H$ of $\HH^*(A)$
satisfying (fg1) and (fg2)?
\end{itemize}
{\em If the answer to (1) is yes, or if one has at hand a version of rank
varieties or other varieties for $A$-modules, one may further ask:}
\begin{itemize}
\item[(2)] Is $V_A(M\ot N)=V_A(M)\cap V_A(N)$
for all finite dimensional $A$-modules $M,N$?
\end{itemize}
{\em The property in (2) above is called the {\em tensor
product property} of varieties for modules. 
The following questions may be asked independently of the first two.}
\begin{itemize}
\item[(3)] For all finite dimensional $A$-modules $M,N$, 
is $M\ot N$ projective if and only if $N\ot M$ is projective?
\item[(4)] For all finite dimensional $A$-modules $M$ and positive
integers $n$,
is $M$ projective if and only if $M^{\ot n}$ is projective?
\end{itemize}
\end{questions}

Note that for a given Hopf algebra $A$, 
if the answers to Questions~\ref{q}(1) and (2) are yes, 
then the answers to (3) and (4) are yes:
By the tensor product property, $V_A(M\ot N)=V_A(N\ot M)$, and by 
Proposition~\ref{prop:properties}, 
this variety has dimension 0 if and only if $M\ot N$ (respectively
$N\ot M$) is projective.
Also by the tensor product property, $V_A(M^{\ot n})=V_A(M)$, and again $M^{\ot n}$ 
(respectively $M$) 
is projective if and only if the dimension of its support variety is 0.

We will see in the next section that there are Hopf algebras for which
the answer to Question~\ref{q}(2) is no, and yet there is another
way to express $V_A(M\ot N)$ in terms of 
$V_A(M)$ and $V_A(N)$. 
So we may wish to consider instead one of the 
following questions about all finite dimensional $A$-modules $M,N$: 

\medskip

\begin{itemize}
\item[(2$'$)] {\em Can $V_A(M\ot N)$ be expressed in terms of $V_A(M)$
and $V_A(N)$?}
\item[or (2$''$)] {\em Is $ \dim V_A(M\ot N) = \dim ( V_A(M)\cap V_A(N))$?}
\end{itemize}

\medskip

\noindent
For either of these questions, if the answer is yes, one may still use support
varieties to obtain valuable information about the tensor product
structure of modules, for example, the property in (2$''$) allows us
to understand the complexity of $M\ot N$ using knowledge of the
support varieties of the tensor factors $M,N$. 

Many mathematicians have worked on Question~\ref{q}(1).
It is closely 
related to a conjecture of Etingof and Ostrik~\cite{Etingof-Ostrik}
that the cohomology ring of a finite tensor category
is finitely generated; this includes the category of finite dimensional
modules of a
finite dimensional Hopf algebra as a special case.
This is condition (fg1$'$).
The further condition (fg2$'$) should follow using similar
proof techniques as for the finite generation of $\coh^*(A,k)$.
One can then take $H$ to be $\coh^{ev}(A,k)\cdot \HH^0(A)$, or
use $\coh^*(A,k)$ directly to define support varieties, as
explained in Section~\ref{sec:var}.
As a cautionary note however, a related conjecture about Hochschild cohomology
of finite dimensional algebras was shown to be false; see, 
e.g.,~\cite{S,Snashall-Solberg,Xu}.

We next discuss some general classes of Hopf algebras for which
the answers to all four Questions~\ref{q} are known to be yes, as well
as those for which some of the four questions are known to have
positive answers, while others remain open. In the next section,
we discuss some classes of Hopf algebras for which the answer
to at least one of the four questions is no.

\subsection*{Finite group algebras}
If $A=kG$, where $G$ is a finite $p$-group and
$\chara(k)=p$, the answers to 
all four Questions~\ref{q} are yes:
As explained in Example~\ref{ex:gpalg}, 
one may take $H=\coh^{ev}(G,k)\cdot \HH^0(kG)$.
Since $\HH^0(kG) \cong Z(kG)$ is a local ring, for the purpose of
defining varieties, this is equivalent to taking $H$ to be
simply $\coh^{ev}(G,k)$, the standard choice (see, e.g.,~\cite{benson91b,Carlson}).
If $G$ is not a $p$-group, the standard 
version of varieties for modules and the one coming from Hochschild
cohomology differ by finite surjective maps. 
The answer to Question~\ref{q}(2)
is yes for the standard version and is no for the Hochschild cohomology version
(tensor a module in a nonprincipal block with the trivial module $k$).
However the answer to the modified question (2$'$) is still yes in this case. 
The answers to Questions~\ref{q}(3) and (4) are yes.

\subsection*{Finite group schemes}
If $A$ is a  finite dimensional 
cocommutative Hopf algebra (equivalently, finite group scheme), 
the answers are also known due to
work of many mathematicians building on work on
finite groups, on restricted Lie algebras~\cite{FP},
and on infinitesimal group schemes~\cite{SFB}.
See, e.g.,~\cite{FP,F-Pevtsova,FS,SFB}, and the surveys~\cite{Friedlander,Pevtsova}.
In this case, one works with support varieties defined via Hopf algebra
cohomology $\coh^*(A,k)$, which is known to satisfy conditions (fg1$'$)
and (fg2$'$)~\cite{FP,FS}, and with
rank varieties, which are homeomorphic to the support varieties. 
The tensor product property is proven using rank varieties~\cite{F-Pevtsova}.
The answers to all four Questions~\ref{q} are yes, under these choices.

\subsection*{Quantum elementary abelian groups}
If $A$ is a  quantum elementary abelian group as defined in 
Example~\ref{ex:QEAG}, the answers to Questions~\ref{q} are known.
Again we work with support varieties defined via Hopf algebra cohomology
$\coh^*(A,k)$ which satisfies (fg1$'$) and (fg2$'$)~\cite{PW1}.
The tensor product property is proven  using rank varieties~\cite{PW}, and 
the answers to  all four Questions~\ref{q} are  yes.

\subsection*{Finite quantum groups and function algebras and more}
If $A$ is a  small quantum group $u_q({\mathfrak{g}})$
(see Example~\ref{ex:sqg}), 
by~\cite{BNPP,ginzburg-kumar93}, the cohomology $\coh^*(A,k)$ is known
to be  finitely generated for most values of the parameters, and 
the answer to Question~\ref{q}(1) is yes.
However Question~\ref{q}(2) is open; it was conjectured by Ostrik~\cite{Ostrik}
who developed support variety theory for these Hopf algebras. 
Likewise for some more general  pointed Hopf algebras with abelian
groups of grouplike elements (see~\cite{MPSW})
and some finite quotients of  quantum function algebras (see~\cite{Gordon}). 
If $A$ is the 12-dimensional  Fomin-Kirillov algebra---a pointed Hopf
algebra with nonabelian group of grouplike elements---then $\coh^*(A,k)$ is 
finitely generated (see~\cite{SV}). 
For all of these important examples, there should be a good support variety
theory, yet  the tensor product property is unknown.

\medskip

Applications of varieties for modules abound, and are well developed
for some of the classes of Hopf algebras described above. 
For example, one can construct modules with prescribed support 
(see~\cite{AI,benson91b,Feldvoss-W,F-Pevtsova}).
Representation type can be seen in the varieties 
(see~\cite{Farnsteiner,Feldvoss-W,Gordon2}).
Some of the structure of the (stable) module category can be understood 
from knowledge of particular subcategories analogous to ideals in a ring,
and these are typically parametrized by support
varieties (see~\cite{BCR,BIK,BIKP2,Benson-Witherspoon,CI,PW}). 

Most of the foregoing discussion focuses on Questions~\ref{q}(1) and (2).
We  now give a general context in which the answer to Question~\ref{q}(4) is known 
to be yes, independently of any variety theory: 
Let $A$ be a finite dimensional Hopf algebra, let $M$ be
an $A$-module for which $M\ot M^*\cong M^*\ot M$,
and let $n$ be a positive integer.
Then $M$ is projective if and only if $M^{\ot n}$ is projective.
This statement is a consequence of rigidity, since rigidity implies that
$M$ is a direct summand of $M\ot M^*\ot M\cong M^*\ot M\ot M$.
See~\cite{Plavnik-Witherspoon} for details. 
Hopf algebras for which the tensor product of modules is commutative
up to isomorphism (such as almost cocommutative or quasitriangular Hopf
algebras) always satisfy this condition, and so 
the answer to Question~\ref{q}(4) is yes for these Hopf algebras.

\section{Some negative answers}

In this section we give examples of Hopf algebras for which the
answer to Question~\ref{q}(1) is yes, and there is a reasonable support variety
theory for which the answer to Question~\ref{q}(2$'$) is yes,  
however the answers to Questions~\ref{q}(2), (2$''$), (3), and (4) are no.

Our first class of examples is from~\cite{Benson-Witherspoon} with Benson.
Let $k$ be a field of  positive characteristic $p$ and let $L$ be a finite $p$-group.
Let $G$ be a finite group acting on $L$ by automorphisms. 
Let
$$
   A = kL\ot k[G] 
$$ 
as an algebra, where $k[G]$ is the linear dual of the group algebra $kG$
as in Example~\ref{ex:dual}.
The comultiplication is not the tensor product comultiplication,
but rather is modified by the group action:
$$
   \Delta(x\ot p_g) = \sum_{\substack{a,b\in G \\ ab=g}}
    (x\ot p_a) \ot ((a^{-1}\cdot x)\ot p_b)
$$
for all $x\in L$, $g\in G$.
The counit and antipode are given by $\varepsilon(x\ot p_g) = \delta_{g,1}$
and $S(x\ot p_g) = (g^{-1}\cdot x^{-1})\ot p_{g^{-1}}$ for all $x\in L$,
$g\in G$. 
This is termed the {\em smash coproduct} of $kL$ and $k[G]$, written 
$A=kL\natural k[G]$. 

Since $\{ 1\ot p_g\mid g\in G\}$ is a set of orthogonal central idempotents in $A$,
any $A$-module $M$ decomposes as a direct sum, 
$$
    M = \bigoplus_{g\in G} M_g , 
$$
where $M_g = (1\ot p_g)\cdot  M$.
Note that each component $M_g$ is itself a $kL$-module by restriction
of action to the subalgebra $kL\cong kL\ot k$ of $A$.
By~\cite[Theorem~2.1]{Benson-Witherspoon}, for any two $A$-modules $M,N$,
\begin{equation}\label{eqn:tp}
   (M\ot N)_g \cong \bigoplus_{\substack{a,b\in G \\ ab=g}}
    M_a \ot ( {}^a N_b) 
\end{equation}
as $kL$-modules, where ${}^a N_b$ is the conjugate $kL$-module that has 
underlying vector space $N_b$ and action 
$x\cdot _a n = (a^{-1}\cdot x)\cdot n $ for all $x\in L$, $n\in N$.

As an algebra, $A$ is a tensor product of $kL$ and $k[G]$, and so its
Hochschild cohomology is 
$$
   \HH^*(A) \cong \HH^*(kL)\ot \HH^*(k[G]) \cong \HH^*(kL)\ot k[G] .
$$
The latter isomorphism occurs since $k[G]$ is semisimple; its Hochschild
cohomology is concentrated in degree~0 where it is isomorphic to the center of the
commutative algebra $k[G]$. Let
$$
    H= \coh^{ev}(L,k) \ot k[G] ,
$$
a subalgebra of $\HH^*(A)$ via the embedding of $\coh^{ev}(L,k)$ into
$\HH^*(kL)$ discussed in Example~\ref{ex:gpalg}.
Then $H$ satisfies (fg1) and (fg2) since $\coh^*(L,k)$ does, 
and we use $H$ to define
support varieties for $A$-modules.
The maximal ideal spectrum of $H$ is 
$$
  {\rm{Max}} (H) \cong {\rm{Max}} ( \coh^{ev}(L,k))\times G .
$$
Define support varieties of  $kL$-modules via $\coh^{ev}(L,k)$
as usual (see Example~\ref{ex:gpalg}), and denote such varieties by $V_L$
in order to distinguish them from varieties for the related $A$-modules.
By~(\ref{eqn:tp}), Proposition~\ref{prop:properties}(ii), and the tensor product property
for finite groups, 
\begin{equation}\label{eqn:tpformula}
   V_A((M\ot N)_g) = \bigcup_{\substack{a,b\in G\\ ab=g}}
   (V_L(M_a) \cap V_L( {}^aN_b) ) \times g  
\end{equation}
for each $g\in G$. 
This formula gives a positive answer to Question~\ref{q}(2$'$), yet it
implies that the answers to Questions~\ref{q}(2), (2$''$), (3),
and (4) are no for some choices of $L$ and $G$:
For example, let $p=2$ and  $L=\Z/2\Z \times \Z/2\Z$,
generated by $g_1$ and $g_2$.
Let $G=\Z/2\Z$, generated by $h$, acting on $L$
by interchanging $g_1$ and $g_2$. 
In this case $\coh^{*}(L,k)\cong k[y_1,y_2]$ with $y_1,y_2$ of degree~1,
and so $V_L(k)$
may be identified with affine space $k^2$. 
Let $U=kL/(g_2-1)$, a $kL$-module.
Then ${}^h U \cong kL/(g_1-1)$. 
Note that $V_L(U)$ may be identified with 
the line $y_1=0$ and $V_L( {}^hU)$ may be identified with the line $y_2=0$.
Now let $M = U\ot kp_h$ and $N=U\ot kp_1$, with $A=kL\ot k[G]$ acting factorwise. 
By the tensor product formula~(\ref{eqn:tpformula}), 
$V_A(N\ot M)$ consists of the line $y_1=0$ paired with the group element $h$,
while $V_A(M\ot N)$ has dimension~0. 
Thus the answers to Questions~\ref{q}(2) and (2$''$) are no. 
By Proposition~\ref{prop:properties}, 
$M\ot N$ is projective while $N\ot M$ is not.
Similarly, $M\ot M$ is projective while $M$ is not. 
More such examples are in~\cite{Benson-Witherspoon}, 
including examples showing that for any positive integer $n$,
it can happen that $M^{\ot n}$
is projective while $M^{\ot (n-1)}$ is not,
and examples of modules $M$ for which $V_A(M^*)\neq V_A(M)$. 
These examples are generalized in~\cite{Plavnik-Witherspoon} with Plavnik to
crossed coproducts $kL\natural_{\sigma}^{\tau} k[G]$ whose
algebra and coalgebra structures are twisted by cocycles $\sigma$, $\tau$.

The above examples are all in positive characteristic.
There are characteristic~0 examples in~\cite{Plavnik-Witherspoon} that
are completely analogous, where the group algebra $kL$ is replaced
by a quantum elementary abelian group as in Example~\ref{ex:QEAG}.
Again we find modules whose tensor product in one order is projective
while in the other order is not, and nonprojective modules with a projective
tensor power. 
In fact, these types of examples are very general, as the following theorem shows.

\begin{thm}\cite{Plavnik-Witherspoon}
Let $A$ be a finite dimensional nonsemisimple Hopf algebra
satisfying (fg1), (fg2), and the tensor product property.
Then $A$ is a subalgebra of a Hopf algebra $K$
satisfying (fg1) and (fg2) 
for which the tensor product
property does not hold.
Moreover, there are nonprojective $K$-modules $M$, $N$ for which
$M\ot M$ and $M\ot N$ are projective, while $N\ot M$ is not projective. 
\end{thm}

One such Hopf algebra is a
smash coproduct $(A\ot A) \natural k^{\Z/2\Z}$ where the nonidentity
element of the group $\Z/2\Z$ interchanges the two tensor factors of $A$.
See~\cite{Plavnik-Witherspoon} for details. 
This is a Hopf algebra for which the answer to Question~\ref{q}(1) is yes,
and so it has a reasonable support variety theory and
Question~\ref{q}(2$'$) has a positive answer. 
However the answers to Questions~\ref{q}(2), (2$''$), (3), and (4) 
are no, just as in our earlier classes of examples in this section.

The positive answers in Section~\ref{sec:questions}
 to Questions~\ref{q}
and the negative answers in this section all point to a larger
question:
What properties of a Hopf algebra ensure positive (respectively, negative)
answers to Questions~\ref{q}?


\begin{thebibliography}{99}

\bibitem{AB} L.\ L.\ Avramov and R.-O.\ Buchweitz,
{\em Support varieties and cohomology over complete intersections},
Invent.\ Math.\ 142 (2000), no.\ 2, 285--318.

\bibitem{AI} L.\ L.\ Avramov and S.\ B.\ Iyengar,
{\em Constructing modules with prescribed cohomological support},
Illinois J.\ Math.\ 51 (2007), no.\ 1, 1--20.

\bibitem{AS82} G.\ S.\ Avrunin and L.\ L.\ Scott,
{\em Quillen stratification for modules},
Invent.\ Math.\ 66 (1982), 277--286. 

\bibitem{BK} B.\ Bakalov and A.\ Kirillov, Jr.,
Lectures on Tensor Categories and Modular Functors,
Vol.\ 21, University Lecture Series, American Mathematical Society, Providence,
RI, 2001.



\bibitem{BNPP} C.\ Bendel, D.\ K.\ Nakano, B.\ J.\ Parshall, and
C.\ Pillen, 
{\em Cohomology for quantum groups via the geometry of the nullcone},
Mem.\ Amer.\ Math.\ Soc.\ 229 (2014), no.\ 1077. 

\bibitem{benson91a} D.\ J.\ Benson, Representations and Cohomology I:
Basic representation theory of finite groups and associative algebras, 
Cambridge University Press, 1991.

\bibitem{benson91b} D.\ J.\ Benson, Representations and Cohomology II:
Cohomology of groups and modules, Cambridge studies in advanced
mathematics  31, Cambridge University Press, 1991.

\bibitem{BCR} D.\ J.\ Benson, J.\ Carlson, and J.\ Rickard,
{\em Thick subcategories of the stable module category},
Fundamenta Mathematicae 153 (1997), 59--80.

\bibitem{BEH} D.\ J.\ Benson, K.\ Erdmann, and M.\ Holloway,
{\em Rank varieties for a class of finite-dimensional local algebras},
J. Pure Appl.\ Algebra 211 (2) (2007), 497--510. 

\bibitem{BG} D.\ J.\ Benson and E.\ L.\ Green,
{\em Non-principal blocks with one simple module},
Q.\ J.\ Math.\ 55 (2004) (1), 1--11. 

\bibitem{BIK} D.\ J.\ Benson, S.\ Iyengar, and H.\ Krause,
{\em Stratifying modular representations of finite groups},
Ann.\ of Math.\ 174 (2011), 1643--1684. 

%\bibitem{BIKP} D.\ J.\ Benson, S.\ Iyengar, H.\ Krause, and J.\ Pevtsova,
%{\em Stratification and $\Pi$-cosupport: Finite Groups},
%arxiv:1505.06628. 

\bibitem{BIKP2} D.\ J.\ Benson, S.\ Iyengar, H.\ Krause, and J.\ Pevtsova,
{\em Stratification for module categories of finite group schemes},
arxiv:1510.06773. 

\bibitem{Benson-Witherspoon} D.\ J.\ Benson and S.\ Witherspoon, 
{\em Examples of support varieties for Hopf algebras with
noncommutative tensor products}, Archiv der Mathematik 102 (2014),
no.\ 6, 513--520. 

\bibitem{BE09} P.\ A.\ Bergh and K.\ Erdmann,
{\em The Avrunin-Scott Theorem for quantum complete intersections},
J.\ Algebra 322 (2009), no.\ 2, 479--488. 

\bibitem{BKN} B.\ D.\ Boe, J.\ R.\ Kujawa, and D.\ K.\ Nakano,
{\em Cohomology and support varieties for Lie superalgebras},
Trans.\ Amer.\ Math.\ Soc.\ 362 (2010), no.\ 12, 6551--6590. 

\bibitem{Carlson} J.\ F.\ Carlson,
{\em The varieties and the cohomology ring of a module},
J.\ Algebra 85 (1983), 104--143.

\bibitem{CI} J.\ F.\ Carlson and S.\ Iyengar,
{\em Thick subcategories of the bounded derived category of a finite group},
Trans.\ Amer.\ Math.\ Soc.\ 367 (2015), no.\ 4, 2703--2717. 

\bibitem{DS} M.\ Duflo and V.\ Serganova, 
{\em On associated variety for Lie superalgebras},
arxiv:0507198. 

\bibitem{EHSST} K.\ Erdmann, M.\ Holloway, N.\ Snashall, \O.\ Solberg,
and R.\ Taillefer,
{\em Support varieties for selfinjective algebras}, 
K-Theory 33 (2004), no.\ 1, 67--87. 

\bibitem{Etingof-Ostrik} P.\ Etingof and V.\ Ostrik, {\em Finite tensor categories},
Mosc.\ Math.\ J.\  4 (2004), no.\ 3, 627--654, 782--783

\bibitem{Ev} L.\ Evens, {\em The cohomology ring of a finite group},
Trans.\ Amer.\ Math.\ Soc.\ 101 (1961), 224--239. 

\bibitem{Farnsteiner} R.\ Farnsteiner,
{\em Tameness and complexity of finite group schemes},
Bull.\ London Math.\ Soc.\ 39 (2007), no.\ 1, 63--70. 

\bibitem{Feldvoss-W} J.\ Feldvoss and S.\ Witherspoon,
{\em Support varieties and representation type of small quantum groups},
Int.\ Math.\ Res.\ Not.\ (2010), no.\ 7, 1346--1362.
{\em Erratum}, Int.\ Math.\ Res.\ Not.\ (2015), no.\ 1, 288-290. 

\bibitem{Friedlander}
E.\ M.\ Friedlander,
{\em Spectrum of group cohomology and support varieties},
J.\ K-Theory 11 (2013), 507--516. 

\bibitem{FP} E.\ M.\ Friedlander and B.\ J.\ Parshall, 
{\em On the cohomology of algebraic and related finite groups},
Invent.\ Math.\ 74 (1983), no.\ 1, 85--117. 

\bibitem{FP86} E.\ M.\ Friedlander and B.\ J.\ Parshall,
{\em Support varieties for restricted Lie algebras},
Invent.\ Math.\ 86 (1986), no.\ 3, 553--562. 

\bibitem{FP05} E.\ M.\ Friedlander and J.\ Pevtsova,
{\em Representation theoretic support spaces for finite group schemes},
Amer.\ J.\ Math.\ 127 (2005), 379--420. Correction: Amer.\ J.\ Math.\ 128 (2006), 1067--1068. 

\bibitem{F-Pevtsova} E.\ M.\ Friedlander and J.\ Pevtsova, 
{\em $\Pi$-supports for modules for finite group schemes over a field},
Duke Math.\ J.\ 139 (2007), no.\ 2, 317--368. 

\bibitem{FS} E.\ M.\ Friedlander and A.\ Suslin, {\em Cohomology of finite
group schemes over a field}, Invent.\ Math.\  127 (1997), no.\ 2, 209--270.

\bibitem{ginzburg-kumar93} V.\ Ginzburg and S.\ Kumar, {\em Cohomology
of quantum groups at roots of unity}, Duke Math.\ J.  69 (1993),
179--198.


\bibitem{golod} E.\ Golod, {\em The cohomology ring of a finite $p$-group},
(Russian) Dokl.\ Akad.\ Nauk SSSR  235 (1959), 703--706.

\bibitem{Gordon} I.\ G.\ Gordon, 
{\em Cohomology of quantized function algebras at roots of unity},
Proc.\ London Math.\ Soc.\ (3) 80 (2000), no.\ 2, 337--359. 

\bibitem{Gordon2} I.\ G.\ Gordon,
{\em Complexity of representations of quantised function algebras and 
representation type}, J.\ Algebra 233 (2000), no.\ 2, 437--482. 

\bibitem{Jantzen} J.\ C.\ Jantzen,
{\em Kohomologie von $p$-Lie-Algebren und nilpotente Elemente},
Abh.\ Math.\ Sem.\ Univ.\ Hamburg 56 (1986), 191--219.

\bibitem{montgomery93} S.\ Montgomery, Hopf Algebras and Their Actions
on Rings, CBMS Conf.\ Math.\ Publ.\  82, Amer.\ Math.\ Soc., 1993.


\bibitem{MPSW} M.\ Mastnak, J.\ Pevtsova, P.\ Schauenburg, and S.\ Witherspoon,
{\em Cohomology of finite dimensional pointed Hopf algebras}, 
Proc.\ London Math.\ Soc.\ (3)  100 (2010), no.\ 2, 377--404.

\bibitem{Ostrik} V.\ V.\ Ostrik,
{\em Support varieties for quantum groups},
Functional Analysis and its Applications 32 (1998), no.\ 4, 237--246.

\bibitem{PS} I.\ Penkov and V.\ Serganova,
{\em The support of an irreducible Lie algebra representation},
J.\ Algebra 209 (1998), no.\ 1, 129--142. 

\bibitem{Pevtsova} J.\ Pevtsova,
{\em Representations and cohomology of finite group schemes}, 
in Advances in Representation Theory of Algebras,
EMS Series of Congress Reports (2013), 231--262. 

\bibitem{PW1} J.\ Pevtsova and S.\ Witherspoon, 
{\em Varieties for modules of quantum elementary abelian groups},
Algebras and Rep.\ Th. {\bf 12} (2009), no.\ 6, 567--595.

\bibitem{PW} J.\ Pevtsova and S.\ Witherspoon,
{\em Tensor ideals and varieties for modules of quantum elementary abelian groups},
Proc.\ Amer.\ Math.\ Soc.\ 143 (2015), no.\ 9, 3727--3741. 



\bibitem{Plavnik-Witherspoon} J.\ Plavnik and S.\ Witherspoon, 
{\em Tensor products and support varieties for some 
noncocommutative Hopf algebras}, arxiv:1611.10285. 


\bibitem{Quillen} D.\ Quillen,
{\em The spectrum of an equivariant cohomology ring: I, II},
Ann.\ Math.\ 94 (1971), 549--572, 573--602.


\bibitem{SchTo} S.\ Scherotzke and M.\ Towers,
{\em Rank varieties for Hopf algebras},
J.\ Pure Appl.\ Algebra 215 (2011), 829--838. 

\bibitem{siegel-witherspoon99} S.\ F.\ Siegel and S.\ J.\ Witherspoon,
{\it The Hochschild cohomology ring of a group algebra}, Proc.\
London Math.\ Soc.\  79 (1999), 131-157.

\bibitem{S} N.\ Snashall, 
{\em Support varieties and the Hochschild cohomology ring modulo
nilpotence},
Proceedings of the 41st Symposium on Ring Theory and Representation
Theory, 68--82, Symp.\ Ring Theory Represent.\ Theory Organ.\ Comm.,
Tsukuba, 2009. 

\bibitem{Snashall-Solberg} N.\ Snashall and \O.\ Solberg,
{\it Support varieties and Hochschild cohomology rings}, Proc.\ London
Math.\ Soc.\ (3)  88 (2004), no.\ 3, 705--732. 

\bibitem{Solberg} \O.\ Solberg,
{\em Support varieties for modules and complexes},
Trends in Representation Theory of Algebras and Related Topics,
Contemp.\ Math.\ 406, 239--270, AMS, Providence, RI, 2006. 

\bibitem{SV} D.\ \c{S}tefan and C.\ Vay,
{\em The cohomology ring of the 12-dimensional Fomin-Kirillov algebra},
Adv.\ Math.\ 291 (2016), 584--620. 


\bibitem{Suarez-Alvarez} M.\ Suarez-Alvarez, {\it The Hilton-Eckmann
argument for the anti-commutativity of cup-products}, Proc.\ Amer.\
Math.\ Soc.\ {\bf 132} (8) (2004), 2241--2246. 

\bibitem{SFB} A.\ Suslin, E.\ M.\ Friedlander, and C.\ P.\ Bendel,
{\em Support varieties for infinitesimal group schemes},
J.\ Amer.\ Math.\ Soc.\ 10 (1997), no.\ 3, 729--759. 



\bibitem{venkov} B.\ B.\ Venkov, {\em Cohomology algebras for some classifying
spaces}, Dokl.\ Akad.\ Nauk.\ SSR  127 (1959), 943--944.


\bibitem{Xu} F.\ Xu, 
{\em Hochschild and ordinary cohomology rings of small categories},
Adv.\ Math.\ 219 (2008), 1872--1893. 


\end{thebibliography}
\end{document}